\documentclass[preprint,12pt]{elsarticle}

\usepackage{natbib}

 \usepackage{graphicx}

\usepackage{amssymb}
\usepackage{amsfonts}
\usepackage{amsmath}
\usepackage{amsthm}
\usepackage{color}

\newtheorem{thm}{Theorem}
\newtheorem{lem}[thm]{Lemma}
\newdefinition{rem}{Remark}
\newtheorem{cor}[thm]{Corollary}

\begin{document}

\begin{frontmatter}

\title{Local polynomial regression\\ based on functional data}

\author[kb]{K. Benhenni}
\ead{Karim.Benhenni@upmf-grenoble.fr}
\author[dd]{D. Degras}
\ead{ddegras@samsi.info}

\journal{Journal of Multivariate Analysis}

\date{\today}

 \address[kb]{Laboratoire LJK UMR CNRS 5224, Universit\'e de Grenoble \\
 UFR SHS, BP. 47\\
 F38040 Grenoble Cedex 09, France}

 \address[dd]{Statistical and Applied Mathematical Sciences Institute\\
 19 T.W. Alexander Drive, P.O. Box 14006\\
Research Triangle Park, NC 27709, USA}


\begin{abstract}
Suppose that $n$ statistical units are observed, each following
the model $Y(x_j)=m(x_j)+ \varepsilon(x_j),\, j=1,...,N,$ where
$m$ is a regression function, $0 \leq x_1 <\cdots<x_N \leq 1$ are
 observation times spaced according to a sampling density $f$, and $\varepsilon$ is a continuous-time
error process having mean zero and regular covariance function.
Considering the local polynomial estimation of $m$ and its
derivatives, we derive asymptotic expressions  for the bias and
variance as $n,N\to\infty$. Such results are particularly relevant
in the context of functional data where essential information is
contained in the derivatives. Based on these results, we deduce
optimal sampling densities, optimal bandwidths and asymptotic
normality of the estimator. Simulations are conducted in order to
compare the performances of local polynomial estimators based on
exact optimal bandwidths, asymptotic optimal bandwidths,  and
cross-validated bandwidths.
\end{abstract}
\end{frontmatter}


\section{Introduction}

Local polynomial smoothing is a popular method for estimating the regression function and its derivatives.
Besides its ease of implementation, this nonparametric method enjoys several desirable statistical properties such as design adaptation, good boundary behavior, and minimax efficiency.
See for instance the monograph of  \cite{fg96} for a thorough introduction to local polynomial methods.
In particular, classical kernel methods like the Nadaraya-Watson estimator or the Gasser-M\"uller estimator are closely connected to local polynomials, as they correspond to local polynomial fitting of order zero.
However, kernel estimators do not share the nice properties of higher order local polynomials listed above.

There is a vast literature on the asymptotics of local polynomial regression estimators under independent measurement errors.
Asymptotic bias and variance expressions can be found in \cite{fg96,rw94,wj95}, among others (see also \cite{gm84} for kernel methods).
Such expressions give important qualitative insights on the large sample properties of estimators.
They also allow to find optimal theoretical bandwidths and devise data-driven methods for the key problem of selecting the bandwidth.
See for instance \cite{fghh96, r97,rsw95}.

In the case of correlated errors, 
\citet{owy01}  give an excellent review of the available asymptotic theory and smoothing parameter selection methods in nonparametric regression.
Local polynomial estimators are studied for instance under mixing conditions in \cite{mf97}, under association in \cite{m03}, and more recently under (stationary) short-range dependent errors in 
\citet{fv01} and \cite{P09}.
Bootstrap and cross-validation methods are developed in \cite{hlp95} to select the bandwidth in the presence of short-range and long-range dependence,
while \cite{fov04} propose a plugin method for short-range dependent errors.

In the functional data setting considered here (that is, when for each statistical unit a whole curve is observed at discrete times), several authors have studied the estimation of a regression function by means of the Gasser-M\"uller kernel estimator.
For instance, in the case of (continuous-time) covariance-stationary error processes, 
\citet{hw86} derive asymptotic bias and variance expansions
and select the bandwidth by optimizing an estimate of the integrated mean squared error based on the empirical autocovariance.
This work is extended to nonstationary error processes with parametric covariance in \cite{fnr97}.
 \citet{br06,br07} derive asymptotic bias and variance expressions when the errors are general nonstationary processes with regular covariance.
In the context of smoothing splines, 
\citet{rs91} propose a cross-validation method for functional data that leaves one curve (instead of one time point) out at a time. The optimality properties  of this method are established in \cite{hw93}.
Asymptotic distributions of local polynomial regression estimators for longitudinal or functional data can be found in \cite{y07}.
\citet{d08} provides consistency conditions for general linear estimators and builds normal simultaneous confidence intervals for the regression function.
Studying the local linear estimation of a univariate or bivariate regression function, \citet{d10} 
elaborates a Central Limit Theorem in the space of continuous functions and applies it  to build simultaneous confidence bands and  tests based on supremum norms.
However, no result seems to be available in the functional data setting for the nonparametric estimation of {\em regression derivatives}.

In this paper we consider the situation where, for each of $n$ statistical units, a curve is observed at the same $N$ sampling points generated by a positive density in some bounded interval , say $[0,1]$.
The data-generating process is the sum of a regression function $m$ and a general error process $\epsilon$.
We are interested in the estimation of $m$ and its derivatives by local polynomial fitting.
The main contributions of this work are as follows. \\
First, under differentiability conditions on the covariance function of $\epsilon$,
we derive asymptotic expressions for the bias and variance of the local polynomial estimator as $n,N\to\infty$.
Note that the bias expansions can be found elsewhere in the literature (e.g. \cite{wj95}) 
as they do not
depend on the stochastic structure of the measurement errors.
The variance expansions, on the other hand, provide new and important convergence results for the estimation of the regression function and its derivatives using functional data.
In particular they highlight the influence of the bandwidth and the covariance structure in the first and second order expansion terms.
Second, we deduce optimal sampling densities
(see e.g. \cite{bc92,c85} 
for other examples of optimal designs)
 as well as optimal bandwidths in a few important cases
(local constant or linear fit of $m$, local linear or quadratic fit of $m'$).
These quantities can be estimated in practice by plugin methods.
Third, we prove, for inference purposes, the asymptotic normality of the estimators.
Fourth, we conduct extensive simulations to compare:
(i) local polynomial fits of different orders, (ii) local polynomial fits based on different bandwidths (exact optimal bandwidth,
asymptotic optimal bandwidth, and cross-validation bandwidth).
The simulations use local polynomial smoothers of order $p=0,1, 2$ to estimate $m$ or $m'$
with different target functions, error processes, and values of $n,N$.
With this numerical study, we try to answer three specific questions:
is there a better order of local polynomial fit to use in a given scenario?
Are the performances of local estimators based on asymptotic optimal bandwidths good enough to justify the development of plug-in methods?
Does the naive approach that consists in using the cross-validated bandwidth to estimate $m^{(\nu)}$ for some $\nu\ge 1$ give reasonable results
(note that, in general, cross-validation aims to produce good bandwidths for the estimation of $m$ and not $m^{(\nu)}$)?

The rest of the paper is organized as follows.
The regression model and local polynomial estimators are defined in Section 2.
The main asymptotic results are contained in Section 3.
The simulation study is displayed in Section 4 and a discussion is provided in Section 5.
Finally, proofs are deferred to the Appendix.


\section{Local polynomial regression}

We consider the statistical problem of estimating a
regression function and its derivatives for a fixed design model. We consider $n$
experimental units, each of them having $N$ measurements of the
response:
\begin{equation}\label{model}
Y_{i}(x_j)=m(x_j)+\varepsilon_{i}(x_j) ,  \qquad  i=1, \ldots,
n ,\quad  j=1, \ldots, N,
\end{equation}
where $m$ is the unknown regression function and the
$\varepsilon_i$ are i.i.d. error processes with mean zero and
autocovariance function $\rho$.

The observation points $x_j ,\, j=1,\ldots,N$ are taken to be
regularly spaced quantiles of a continuous positive density $f$ on
[0,1]:
\begin{equation}\label{sampling design}
  \int^{x_j}_0 f(x) dx
= \frac{j-1}{N-1}, \quad
   j = 1, \ldots, N.
\end{equation}
Note that the uniform density $f = 1_{\left[0,1\right]}$ corresponds to  an equidistant design.

Let ${\bar Y}_{\cdot j}=\frac{1}{n}\sum_{i=1}^n Y_i(x_j)$ be the
sample average  at location $x_j$ and let $0 \le \nu \le p$ be
integers. For each $x\in [0,1]$, the local polynomial estimator of
order $p$ of the existing $\nu$th order derivatives $m^{(\nu)}(x) $ of the regression function 
 is defined as $\hat{m}_\nu(x) =\nu! \,
\widehat{\beta}_\nu(x)$, where $\widehat{\beta}_\nu(x)$ is the
$\nu$th component of the estimate
$\boldsymbol{\widehat{\beta}}(x)=(\widehat{\beta}_0(x),\ldots,\widehat{\beta}_p(x))$
 which is the solution to the minimization problem\begin{equation*}
\min_{\boldsymbol{\beta}(x)} \:\:\displaystyle \sum_{j=1}^N
\left({\bar Y}_{\cdot j}- \sum_{k=0}^p \beta_k(x)(x_j-x)^k
\right)^2 \frac{1}{h} \,K\left( \frac{x_j-x}{h}\right)  .
\end{equation*}

where $h$ denotes a positive bandwidth and $K$ is a kernel
function. Let $\bar{\mathbf{Y}}= (\bar{Y}_{\cdot 1},\ldots,
\bar{Y}_{\cdot N})'$ and denote the canonical basis of $
\mathbb{R}^{p+1} $ by
 $(\mathbf{e}_k)_{ k=0,\ldots,p}$
($\mathbf{e}_k$ has a 1 in the $(k+1)$th position and
 0 elsewhere). Finally define the matrix
$$
\mathbf{X}=
\left(\begin{array}{cccc}
1 &  (x_1-x) & \cdots & (x_1-x)^q \\
\vdots & \vdots &   & \vdots \\
1 &  (x_N-x) & \cdots & (x_N-x)^q
\end{array}
\right)
$$
and $\mathbf{W}=\mathrm{diag}\left( \frac{1}{h}\,
K\left(\frac{x_j-x}{h}\right)\right) $. Then the estimator
$\hat{m}_\nu(x)$ can be written  as
\begin{equation}\label{def estimateur}
\hat{m}_\nu(x)  =  \nu!\, \mathbf{e}_\nu' \, \widehat{\boldsymbol{\beta}}(x) ,\quad \textrm{with }
 \widehat{\boldsymbol{\beta}}(x)= \mathbf{(X'W X)}^{-1}  \mathbf{X'W}\bar{\mathbf{Y}}.
\end{equation}


\section{Asymptotic study}

\subsection{Bias and variance expansions}

The following assumptions are needed for 
the asymptotic study of  $\hat{m}_\nu(x)$:

\begin{enumerate}
\item[(A1)] The kernel $K$ is a Lipschitz-continuous, symmetric
 density function with support $[-1,1]$.
\item[(A2)] The bandwidth $h=h(\nu,n,N)$ satisfies $h\to 0 $,
$Nh^2\to \infty$, and $nh^{2\nu}\to\infty$ as $n,N\to\infty$.
\item[(A3)] The regression function $m$ has $(p+2)$ continuous
derivatives on $\left[0,1\right]$. \item[(A4)] The sampling
density $f$ has one continuous derivative on $\left[0,1\right]$.
    \item[(A5)]The covariance function $\rho$ is continuous on the unit square $[0,1]^2$ and
has continuous first-order partial derivatives off the main diagonal.
These derivatives have left and right limits on the main diagonal determined by
    $\rho^{(0,1)}(x,x^-)=\lim_{y\nearrow x} \rho^{(0,1)}(x,y)
    \hbox{ and }\rho^{(0,1)}(x,x^+)=\lim_{y\searrow x}  \rho^{(0,1)}(x,y) . $
\end{enumerate}

We now introduce several useful quantities associated to $K$. 
Let  $\mu_k= \int_{-1}^1 u^k K(u)du$ be the $k$th moment
of $K$ and the vectors $\mathbf{c}= (\mu_{p+1},\ldots,
\mu_{2p+1})'$ and $ \mathbf{\tilde{c}}= (\mu_{p+2},\ldots,
\mu_{2p+2})'$. Let  $\mathbf{S}= (\mu_{k+l})$,
$\mathbf{\tilde{S}}= (\mu_{k+l+1})$, $\mathbf{S^{\ast}}=(\mu_k
\mu_l)$, and $\mathbf{A}=\left( \frac{1}{2} \iint_{[-1,1]^2}
|u-v| u^k v^l K(u) K(v)du dv  \right)$ be matrices of size
$(p+1)\times (p+1)$ whose  elements are indexed by $k,l=0,\ldots,p$.

\bigskip 

The asymptotic bias and variance of the estimator $\hat{m}_\nu(x)$,
 for a given $x\in(0,1)$ and $\nu\in \{ 0,\ldots,p \} $,
are established in the following theorem.

\begin{thm}\label{t1}
Assume (A1)--(A5). Then as $n,N \to \infty$,
\begin{align*}
 \mathrm{Bias}(\hat{m}_\nu(x)) & =  \frac{\nu!\,m^{(p+1)}(x)}{(p+1)!} \left(
\mathbf{e}_{\nu}' \mathbf{S^{-1}}  \mathbf{c} \right) h^{p+1-\nu} + o(h^{p+2-\nu})
 \\
\noalign{ $\displaystyle + \nu!  \left\{  \frac{m^{(p+2)}(x)}{(p+2)!}  \mathbf{e}_{\nu}' \mathbf{S}^{-1}\mathbf{\tilde{c}}
+ \frac{m^{(p+1)}(x)}{(p+1)!}\frac{f'(x)}{f(x)}
\left(\mathbf{e}_{\nu}' \mathbf{S}^{-1}  \mathbf{\tilde{c}}
 -\mathbf{e}_{\nu}' \mathbf{S}^{-1} \mathbf{\tilde{S}} \mathbf{S}^{-1}  \mathbf{c} \right)\right\}h^{p+2-\nu} $ \vspace*{1mm}} 
\noalign{\noindent  and \vspace*{2mm}} 
\mathrm{Var}(\hat{m}_\nu (x))
& = \frac{(\nu!)^2 \rho(x,x)}{nh^{2\nu}}\: \mathbf{e}_{\nu}' \mathbf{S}^{-1} \mathbf{S}^\ast  \mathbf{S}^{-1} \mathbf{e}_{\nu}+ o\left(\frac{1}{nh^{2\nu-1}}\right)  \\
& \quad - \frac{(\nu!)^2}{nh^{2\nu-1}} \left( \rho^{(0,1)}(x,x^-)
- \rho^{(0,1)}(x,x^+)   \right)
  \mathbf{e}_{\nu}' \mathbf{S}^{-1}\mathbf{AS}^{-1} \mathbf{e}_{\nu}  .
\end{align*}
\end{thm}

\begin{rem}\label{th1 and literature}
The bias  expansion of Theorem \ref{t1}
 does not depend on the nature of  the measurement errors (continuous time processes).
A similar expansion can be found e.g. in \cite{fg96} 
in the context of independent errors.
Also, the present variance expansion extends the
results of \cite{br07,hw86} 
on the nonparametric estimation of the regression function with
the Gasser-M\"uller estimator.
\end{rem}

\begin{rem}\label{rem:2nd order}
The reason for presenting second-order expansions in Theorem
\ref{t1} is that first-order terms may vanish due to the symmetry
of $K$ which causes its odd moments to be null. For
instance, the first-order terms  in the bias and the variance
vanish, respectively, whenever $p-\nu$ is even and $\nu$ is odd. In
both cases, the second-order terms generally allow to find exact
rates of convergence and asymptotic optimal bandwidths as in
Corollary \ref{local opt h}.
\end{rem}

If the covariance $\rho$ has continuous first
derivatives at $(x,x)$, the second-order variance term in Theorem
\ref{t1} vanishes since $\rho^{(0,1)}(x,x^-)=\rho^{(0,1)}(x,x^+)$. Thus, the variance expansion does not depend
on $h$ when $\nu=0$ or $\nu$ is odd (see Remark \ref{rem:2nd
order}). This makes it impossible to assess the effect of
smoothing on the variance of $\hat{m}_\nu(x)$ nor to optimize the
mean squared error with respect to $h$. This problem can be solved
by deriving higher-order variance expansions under stronger
differentiability assumptions on $f$ and $\rho$. As general
higher-order expansions are quite messy and difficult to
interpret, we restrict ourselves to the central case of an equidistant sampling design with $f = 1_{[0,1]}$. For this
purpose, we introduce the matrices $\mathbf{A}_1 =
\left(\frac{1}{2}(\mu_k\mu_{l+2} + \mu_{k+2}\mu_l) \right) $, $
\mathbf{A}_2 =  (\mu_{k+1} \mu_{l+1})$,
 and $  \mathbf{A}_3 =  \big(\frac{1}{6}(\mu_{k+3} \mu_{l+1}  + \mu_{k+1}\mu_{l+3})\big)$ indexed by $k,l=0,\ldots,p$.

\begin{thm}\label{t2}
Assume (A1)--(A5) with $f\equiv1$ on $[0,1] $ (equidistant design).

\begin{itemize}
\item Case $\nu$ even. Assume further that $\rho $ is twice
continuously differentiable at $(x,x) $ and $Nh^3
\to \infty$
 as $n,N \to \infty$. Then
\begin{equation*}
\begin{split}
\mathrm{Var}  ( \hat{m}_{\nu}(x)) &  = \frac{(\nu!)^2 \rho(x,x)}{nh^{2\nu}} \mathbf{e}_\nu'  \mathbf{S}^{-1}\mathbf{S}^\ast \mathbf{S}^{-1}   \mathbf{e}_\nu  \\
& \qquad + \frac{(\nu!)^2\rho^{(0,2)}(x,x)}{nh^{2\nu -2}} \mathbf{e}_\nu'  \mathbf{S}^{-1} \mathbf{A}_1 \mathbf{S}^{-1}   \mathbf{e}_\nu 
   + o\left( \frac{1}{nh^{2\nu -2}} \right) .
\end{split}
\end{equation*}

\item Case $\nu$ odd. Assume further that $\rho $ is four times continuously differentiable at  $(x,x)$ and $Nh^5 \to \infty$
 as $n,N \to \infty$. Then
\begin{equation*}
\begin{split}
\mathrm{Var}  ( \hat{m}_{\nu}(x)) &  =  \frac{(\nu!)^2 \rho^{(1,1)}(x,x)}{nh^{2\nu -2}} \mathbf{e}_\nu'  \mathbf{S}^{-1} \mathbf{A}_2 \mathbf{S}^{-1}   \mathbf{e}_\nu
\\ 
& \qquad +   \frac{(\nu!)^2 \rho^{(1,3)}(x,x)}{nh^{2\nu -4}} \mathbf{e}_\nu'  \mathbf{S}^{-1} \mathbf{A}_3 \mathbf{S}^{-1}   \mathbf{e}_\nu   + o\left( \frac{1}{nh^{2\nu -4}} \right).
 \end{split}
\end{equation*}
\end{itemize}
\end{thm}

\begin{rem}
The function $m^{(\nu)}$ can be estimated consistently without smoothing the data in model \eqref{model}.
Interpolation methods would also be consistent, as can be checked from \citet{d08}. 
Moreover, looking at Theorems \ref{t1} and \ref{t2}, it is not
clear whether the variance of $\hat{m}_\nu (x)$ is a decreasing
function of $h$. In other words, smoothing more may not always
reduce the variance of the estimator. 
See \citet{car00} for a 
similar observation in the context of functional principal
components analysis.
\end{rem}


\subsection{Optimal sampling densities and bandwidths}
\label{optimal f and h}

In this section, we discuss the optimization of the (asymptotic) mean squared error
\begin{align*}
\mathrm{MSE} & =\mathbb{E} (\hat{m}_\nu(x)-m^{(\nu)}(x) )^2\\
&=  \mathrm{Bias}(\hat{m}_\nu(x))^2 + \mathrm{Var}(\hat{m}_\nu (x))
\end{align*}
in Theorem \ref{t1} with respect to the sampling density $f$ and
the bandwidth $h$. A similar optimization could be carried out in
 Theorem \ref{t2} where the covariance function $\rho$ is
assumed to be more regular (twice or four times differentiable).

\smallskip

We first examine the choice of $f$ that minimizes the asymptotic
squared bias of $\hat {m}_\nu(x)$ since that the asymptotic
variance of $\hat {m}_\nu(x)$ is independent of $f$, as can be
seen in Theorem \ref{t1}. This optimization may be useful in
practice especially when the grid size $N$ is not too large and
subject to a sampling cost constraint.

For $p-\nu$ even, $\mathbf{e}_{\nu}' \mathbf{S^{-1}} \mathbf{c}=0$
so that the first-order term in the bias vanishes, as noted in
Remark \ref{rem:2nd order}. Moreover, the second-order term can be
rendered equal to zero (except at zeros of $m^{(p+1)}(x)$) by
taking a sampling density $f$ such that
 $g_{p,\nu}(x)=0$, where
\begin{equation*}
g_{p,\nu}(x) = \frac{m^{(p+2)}(x)}{(p+2)!}  \mathbf{e}_{\nu}' \mathbf{S}^{-1}\mathbf{\tilde{c}}
+ \frac{m^{(p+1)}(x)}{(p+1)!}\frac{f'(x)}{f(x)}
\left(\mathbf{e}_{\nu}' \mathbf{S}^{-1}  \mathbf{\tilde{c}}
 -\mathbf{e}_{\nu}' \mathbf{S}^{-1} \mathbf{\tilde{S}} \mathbf{S}^{-1}  \mathbf{c} \right) .
\end{equation*}
The solution of the previous equation is
\begin{equation}\label{f optim}
 f_0(x) = d_0^{-1}   \left| m^{(p+1)}(x) \right|^{\gamma/(p+2)} ,
\end{equation}
with $d_0$ such that $\int_0^1 f_0(x)dx = 1$ and $\gamma = \frac {\mathbf{e}_{\nu}' \mathbf{S}^{-1}\mathbf{\tilde{c}}}{ \left(
 \mathbf{e}_{\nu}' \mathbf{S}^{-1} \mathbf{\tilde{S}} \mathbf{S}^{-1}  \mathbf{c} -\mathbf{e}_{\nu}' \mathbf{S}^{-1}  \mathbf{\tilde{c}} \right)} $.
Observe that $f_0(x)$ is well-defined over $[0,1]$ if and only if  $m^{(p+1)}(x) \ne 0$ for all $x\in[0,1]$.

With the choice $f=f_0$, the bias of $\hat{m}_\nu(x)$ is of order $ o(h^{p+2-\nu})$, so that a higher order expansion would be required to get the exact rate of convergence.
In practice, the density $f_0$ depends on the unknown quantity $m^{(p+1)}(x)$. However, an approximation of $f_0(x)$ can be obtained by replacing in \eqref{f optim}  the derivative $m^{(p+1)}(x)$ by a local polynomial estimator $\hat{m}_{(p+1)}(x)$.

For $p-\nu$ odd, the first-order term in the bias is non zero (if
$m^{(p+1)}(x)\ne 0 $) but does not depend on $f$. On the other
hand, the second-order term vanishes for any sampling density $f(x)$.
Therefore, a higher order expansion of the bias would be required
to get exact terms that depend on $f(x)$ and could then be
optimized.

\medskip

We turn to the optimization of the bandwidth $h$ and start with a
useful lemma whose proof is in the Appendix.
\begin{lem}\label{alpha nonnegative}
Assume (A5) and define   $\alpha(x) =  \rho^{(0,1)}(x,x^-) - \rho^{(0,1)}(x,x^+)$. Then $\alpha(x) \ge 0$.
\end{lem}
This lemma is easily checked for covariance-stationary processes
(see e.g. \citet{hw86}) 
but is less intuitive for general covariance functions $\rho$. It can be helpful in
determining whether the asymptotic variance of $\hat{m}_\nu(x)$ is
a decreasing function of $h$, in which case the MSE can be
optimized. More precisely, in order to derive asymptotic optimal
bandwidths throughout this section, we need to assume that
$\alpha(x)>0$ (or assume higher order differentiability for $\rho$  if $\alpha(x)=0$;  
see Remark \ref{optimization for more regular covariance}).

When estimating the regression function itself, $\nu=0$, the
leading variance term in Theorem \ref{t1} does not depend on $h$.
If $\mathbf{e}_0' \mathbf{S}^{-1}\mathbf{AS}^{-1}\mathbf{e}_0 \le
0 $, then the second-order variance term (in $h/n$) is nonnegative
and the optimization of the MSE yields the solution $h=0$, which
is not admissible. In fact, we suspect that $ \mathbf{e}_0'
\mathbf{S}^{-1}\mathbf{AS}^{-1}\mathbf{e}_0
>0 $ for all kernels $K$ satisfying (A1) and all integers $p\ge 0$,
although we have only checked it for the special and interesting
cases $p \le 2$ (local constant, linear, or quadratic fit) but
cannot provide a proof for more general $p$. Proceeding with this
conjecture, the optimal bandwidth for the MSE exists and can be
obtained from Theorem \ref{t1}. However, some caution must be
taken to separate the cases $p$ even and $p$ odd, for which the
bias expressions are different (see Remark \ref{rem:2nd order} and
the optimization of $f$ above). More precisely,  if $\nu =0$ and
$p$ is odd, then the asymptotic optimal bandwidth is
$$
h_{opt}=\left\{ \frac{ (p+1)!^2 \, (  \mathbf{e}_0' \mathbf{S}^{-1}\mathbf{AS}^{-1}\mathbf{e}_0)\,  \alpha(x)
}{ (2p+2) \left( m^{(p+1)}(x)\right)^2  \left( \mathbf{e}_{0}' \mathbf{S^{-1}}  \mathbf{c}\right)^2}\right\}^{1/(2p+1)} n^{-1/(2p+1)}.
$$

In the case where $\nu=0$ and $p$ is even, the asymptotic optimal bandwidth becomes
$$
h_{opt}=\left\{ \frac{(  \mathbf{e}_0' \mathbf{S}^{-1}\mathbf{AS}^{-1}\mathbf{e}_0)\, \alpha(x) }{(2p+4)\, g_{p,0}(x)^2  }\right\}^{1/(2p+3)} n^{-1/(2p+3)} .
$$
Note that in the above optimization, it is assumed that $g_{p,0}(x)\ne 0$, that is, $f$ is different of the optimal sampling density $f_0$.
To optimize the MSE when the optimal density $f_0$ is used, then
it would be necessary in this case to derive a higher order
expansion for the bias, and it can be shown that the optimal
bandwidth would then be of order at least $n^{-1/(2p+5)}$.

In the following corollary, we give the optimal bandwidth
$h$ in two important cases ($\nu\in \{0,1\}$), using, for simplicity, a uniform
sampling density $f\equiv 1$ on $[0,1]$. 
Optimal bandwidths can be obtained similarly for $\nu \ge 2$.

\begin{cor}\label{local opt h}
Assume  (A1)-(A5) with $f\equiv 1$ on $[0,1]$ and $\alpha(x)>0$.
\begin{enumerate}
\item[1.] Local constant or linear estimation of $m$ ($\nu=0$, $p\in \{0,1\}$). Assume further that $ m''(x) \ne 0$
and $Nn^{-2/3} \to \infty$ as $n,N\to\infty$. Then the optimal bandwidth for the asymptotic MSE of $\hat{m}_0 (x)$ is
\begin{equation*}
h_{opt} = \left(\frac{ \alpha(x)}{2\mu_2^{2}\, m''(x)^{2}} \iint_{\mathbb{R}^2} |u-v| K(u)K(v)dudv \right)^{1/3}    n^{-1/3}.
\end{equation*}
\item[2.] Local linear or quadratic estimation of $m'$ $(\nu=1$, $p\in\{1,2\}$). Assume further that $m^{(3)}(x)\ne 0$
and $Nn^{-2/5} \to \infty$ as $n,N\to\infty$. Then the optimal bandwidth for the asymptotic MSE of $\hat{m}_1 (x)$ is
\begin{equation*}
h_{opt} =\left( - \frac{9  \, \alpha(x) }{ 2\mu_4^2 \, m^{(3)}(x)^2} \iint_{\mathbb{R}^2} |u-v| uv K(u)K(v)dudv  \right)^{1/5} n^{-1/5} .
\end{equation*}
\end{enumerate}
\end{cor}

Corollary \ref{local opt h} provides the theoretical basis for a
plug-in method to select $h$. Developing such a method and
studying its theoretical properties is however beyond the scope of
this paper. In Section \ref{sec: num study}, the optimal
bandwidths of Corollary  \ref{local opt h} are used as benchmarks
to assess the popular cross-validation procedure.

\begin{rem}
In the cases ($\nu=0,p=1$) and ($\nu=1,p=2$) of Corollary \ref{local opt h},
the results actually hold for any sampling density $f$ satisfying (A4).
Also, the first part of the corollary corresponds to Theorem 3 of \cite{hw86}
and Corollary 2.1 of \cite{br07} when the Gasser-Muller
estimator is used along with an equidistant sampling.
\end{rem}

\begin{rem}\label{optimization for more regular covariance}
Under the assumptions of Theorem \ref{t2}, in case 1 of Corollary
\ref{local opt h}, the optimal bandwidth is $h_{opt} = \left(-
\frac{2  \rho^{(0,2)}(x,x)  }{\mu_2 m''(x)^2}\right)^{1/2}
n^{-1/2}$ if $\rho^{(0,2)}(x,x)<0$ and $Nn^{-3/2}\to\infty$ as
$n,N\to\infty$, otherwise the optimization will not be possible.
 Likewise, in case 2, the optimal bandwidth is $h_{opt} =  \left(-  \frac{6 \mu_2 \rho^{(1,3)}(x,x)  }{\mu_4 m^{(3)}(x)^2}\right)^{1/2} n^{-1/2} $
 provided that $\rho^{(1,3)}(x,x)<0$ and $Nn^{-5/2}\to\infty$ as $n,N\to\infty$.
\end{rem}

\begin{rem}\label{optimization IMSE}
Theorem 1 can also be harnessed  to derive optimal bandwidths for  global error measures such as the integrated mean squared error
\begin{equation}\label{IMSE}
\mathrm{IMSE} = \int^1_0 \mathbb{E} \left(\hat{m}_\nu(x)-m^{(\nu)}(x)\right)^2 w(x)dx
\end{equation}
where $w$ is a bounded, positive weight function. More precisely, denoting by $[-\tau,\tau]$ the support of $K$, the bias and variance
expansions in Theorem 1 hold uniformly over $ [\tau h,1-\tau h]$, and their convergence rates are maintained in the boundary regions $[0,\tau h)$ and $(1-\tau h,1]$ (only the multiplicative constants are lost). 
As $n,N\to \infty$, the IMSE is therefore equivalent to the weighted integral over $[0,1]$ of the (squared) bias plus variance expansions of Theorem 1.
One can thus replace the terms $\alpha(x)$ and $(m^{(\nu)}(x))^2$  in Corollary \ref{local opt h} by $\int_0^1 \alpha(x)w(x)dx$ and $\int_0^1 (m^{(\nu)}(x))^2 w(x)dx$, respectively, to obtain global optimal bandwidths.
\end{rem}

\begin{rem}\label{estimation int_alpha_x}
For the implementation of the (global) optimal bandwidths corresponding to Corollary \ref{local opt h},
the integral $\int^{1}_{0}\alpha(x)w(x) dx$ can be estimated by  
$ V_N = \frac{1}{n} \sum_{i=1}^n \sum_{j=2}^N \left(Y_i(x_j)-Y_i(x_{j-1})\right)^2  w(x_j)$, 
namely the quadratic variations of the sample processes $Y_i$.
The previous estimator is almost surely consistent as $n,N\to\infty$; see e.g. \cite{P99} in the case of Gaussian processes.
\end{rem}


\subsection{Asymptotic normality}

The asymptotic bias and variance expansions of Theorems \ref{t1} and \ref{t2}
provide the centering and scaling
required to determine the limit distribution of the estimator $ \hat{m}_\nu(x)$.
Besides, one may observe that 
$ \hat{m}_\nu(x)= \frac{1}{n} \sum_{i=1}^n \hat{m}_{\nu,i}(x)$,
where the $\hat{m}_{\nu,i}$ are the local polynomial smoothers of
the curves $Y_i, \, i=1,\ldots, n$. Since the
$\hat{m}_{\nu,i}(x)$'s are i.i.d. with finite variance as the
$Y_i$'s, the Central Limit Theorem can then be applied to
$\hat{m}_{\nu}(x)$
 as $n,N\to\infty$. We consider the asymptotic distribution of
 $\hat{m}_{\nu}(x)$ according to the parity
 of $\nu$ (see Remark \ref{rem:2nd order} on the vanishing terms in the
 asymptotic variance)
in order to get the correct scaling term,
 and also to ensure that the bias term is asymptotically negligible
 when multiplied by the scaling rate. The latter is guaranteed by imposing
 extra conditions on the bandwidth $h$.
 Denoting the convergence in distribution
 by $\stackrel{d}{\to}$ and  the centered normal distribution with variance $\sigma^2$
 by $\mathcal{N}(0,\sigma^2)$, the limit distribution of $ \hat{m}_\nu(x)$ is given by the following result.
\begin{thm}\label{cor: as-norm}  Assume (A1)-(A5).
\begin{itemize}
\item Case $\nu$ even. Assume further that $nh^{2p+4}\to 0$ if $p$ is even, resp.  $nh^{2p+2}\to 0$ if $p$ is odd, as $n,N \to \infty$.
Then
\[  \sqrt{n  h^{2\nu}} ( \hat{m}_\nu(x) - m^{(\nu)}(x))  \stackrel{d}{\to} \mathcal{N} \left( 0, (\nu!)^2 \rho(x,x)
  ( \mathbf{e}_{\nu}' \mathbf{S}^{-1} \mathbf{S}^\ast  \mathbf{S}^{-1} \mathbf{e}_{\nu})\right) .  \]
\item  Case $\nu$ odd and $\alpha(x)>0$. Assume further that $nh^{2p+1}\to 0$ if $p$ is even, resp.  $nh^{2p+3}\to 0$ if $p$ is odd, as $n,N \to \infty$.
Then
\begin{equation*}
\sqrt{ n h^{2\nu-1}}   ( \hat{m}_\nu(x) - m^{(\nu)}(x) )  \stackrel{d}{\to} \mathcal{N} \left( 0,
  (\nu!)^2 \alpha(x) |\mathbf{e}_{\nu}' \mathbf{S}^{-1}
\mathbf{A}    \mathbf{S}^{-1} \mathbf{e}_{\nu} |  \right)  .
\end{equation*}
\item  Case $\nu$ odd and $\alpha(x)=0$ (i.e. $\rho$ is continuously differentiable in a neighborhood of $(x,x)$). 
Assume further that $\rho$ is four times differentiable at $(x,x)$,   
 that $Nh^5\to\infty$, 
and that $nh^{2p}\to 0$ if $p$ is even, resp.  $nh^{2p+2}\to 0$ if $p$ is odd, as $n,N \to \infty$.
Then
\begin{equation*}
\sqrt{ n h^{2\nu-2}}   ( \hat{m}_\nu(x) - m^{(\nu)}(x) )  \stackrel{d}{\to} \mathcal{N} \left( 0,
  (\nu!)^2 \rho^{(1,1)}(x,x) (\mathbf{e}_\nu'  \mathbf{S}^{-1} \mathbf{A}_2 \mathbf{S}^{-1}   \mathbf{e}_\nu)  \right)  .
\end{equation*}
\end{itemize}
\end{thm}


\section{Numerical study}\label{sec: num study}

In this section we compare the numerical performances of local
polynomial estimators based on several global bandwidths.
Specifically, we examine the local constant and linear estimation
of a regression function $m$ and the local linear and quadratic
estimation of $m'$. We consider three types of bandwidths: (i) the global
versions of the asymptotic optimal bandwidths of Corollary
\ref{local opt h}, denoted by $h_{as}$; (ii) the bandwidths that
minimize the IMSE on a finite sample with the target $m$ or $m'$
and weight function $w\equiv 1$ in \eqref{IMSE}, which are called 
exact optimal bandwidths and denoted by $h_{ex}$; (iii) the popular
``leave-one-curve-out" cross-validation bandwidths, denoted by
$h_{cv}$. The interested reader may refer to Rice and Silverman
\cite{rs91} and Hart and Wehrly \cite{hw93} for a detailed account
and theoretical justification of this cross-validation method. In
short, this method selects the smoothing parameter for which the
estimator based on all observed curves but one
 predicts best the remaining curve.
In  model \eqref{model}, considering
the  local polynomial estimation of $m$ by $\hat{m}_{0}$ in \eqref{def estimateur},
$h_{cv}$ is obtained by minimizing  
\begin{equation}\label{CV}
\mathrm{CV}(h)= \frac{1}{nN} \sum_{i=1}^n \sum_{j=1}^N \left( \hat{m}_{0}^{(-i)}(x_j;h) - Y_{i}(x_j)\right)^2
\end{equation}
where $\hat{m}_{0}^{(-i)}(\cdot;h)$ is the local polynomial smoother of order $p$ and bandwidth $h$ applied to the data
$Y_{k}(x_j), \, k \in \{1,\ldots,n \} \setminus \{ i\} , \, j=1,\ldots,N$.
Although the cross-validation score \eqref{CV} is designed for the estimation of $m$,
it is of interest to see how it performs for the estimation of derivatives $m^{(\nu)}, \nu \ge 1$.
This question is particularly justified in local polynomial fitting,
where all derivatives of $m$ up to order $p$ are being estimated simultaneously.

\begin{figure}[d]
\begin{center}
\includegraphics[scale=.75]{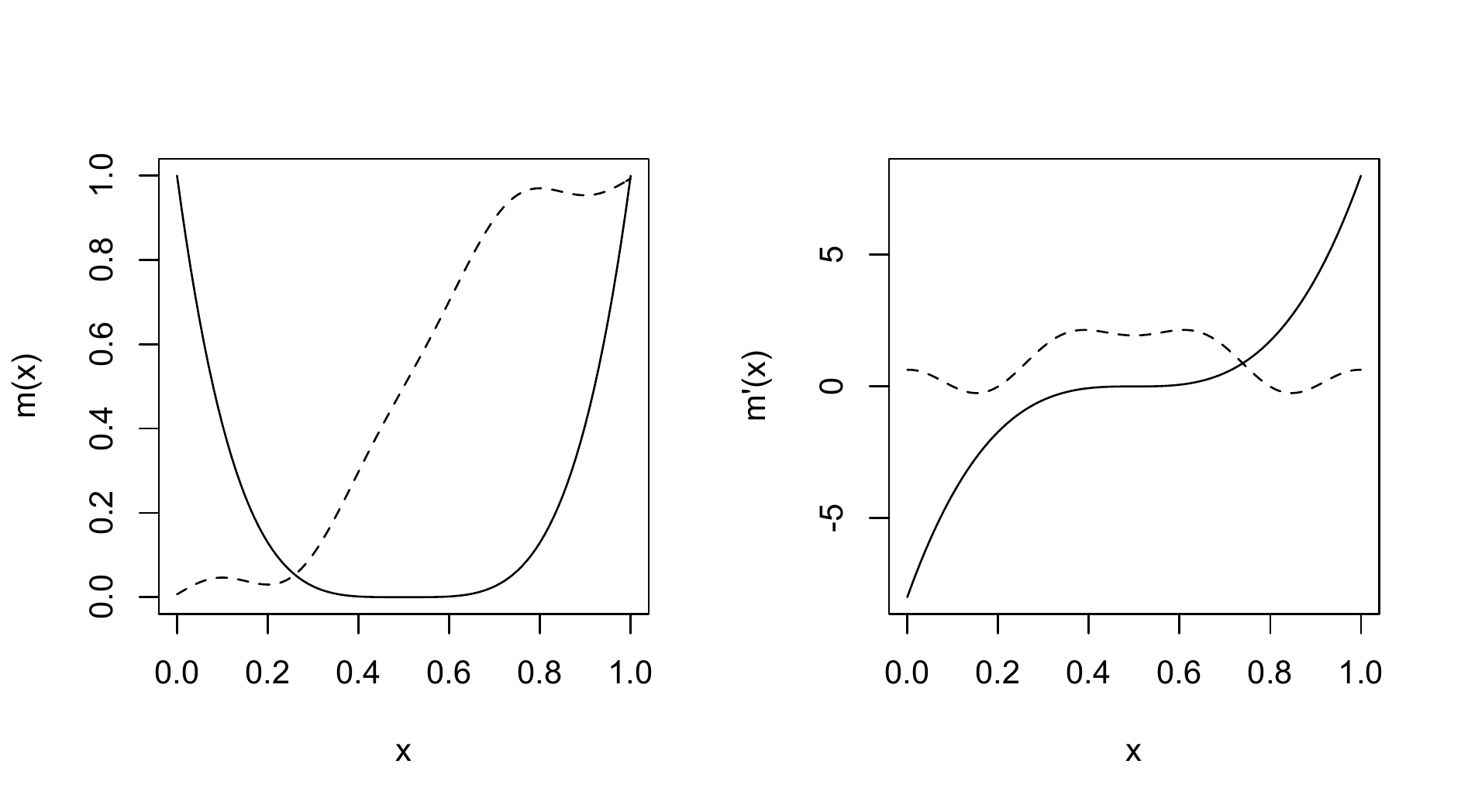}
\caption{Left panel: regression functions $m_1$ (solid line) and $m_2$ (dashed line). Right panel: first derivatives $m_1'$ (solid) and $m_2'$ (dashed).}
\label{m1m2}
\end{center}
\end{figure}

The regression functions used in the simulations are
\begin{equation}
\left\{
\begin{array}{l}
\displaystyle m_1(x)  = 16 (x-0.5)^4, \\
\displaystyle  m_2(x)  = \frac{1}{1+e^{-10(x-0.5)}} + 0.03 \sin(6\pi x).
\end{array}
\right.
\end{equation}
The polynomial function $m_1$ has unit range and has relatively
high curvature away from its minimum at $x=0.5$.
The function $m_2$ is a linear combination of a logistic function
and a rapidly varying sine function. The factor $0.03$ is chosen
so that the sine function has small influence on $m_2$ but a much
larger on $m_2'$. These functions and their first derivatives are
displayed in Figure \ref{m1m2}.

\bigskip

\begin{center}
$>>>>>>$ \qquad \quad INSERT FIGURE 1 ABOUT HERE  \quad \qquad $<<<<<<$
\end{center}

\bigskip

For the stochastic part of \eqref{model} we consider Gaussian processes with mean zero and covariance functions
\begin{equation}
\left\{
\begin{array}{l}
\rho_1(x,y)  = \min (x,y), \\
\rho_2(x,y)  = e^{-15 |x-y|}\, .
\end{array}
\right.
\end{equation}

\noindent The first error process is a standard Wiener process on $[0,1]$; the second is a stationary Ornstein-Uhlenbeck process.
 The parameter $\lambda=15$ in $\rho_2$ allows to inspect various correlation levels
 between two consecutive measurements (e.g. 0.22 for $N=10$ and 0.86 for $N=100$).
The variance levels of these processes are chosen so that the signal-to-noise ratios
($\mathrm{SNR}= \left\{ \max_t m(t) -\min_t m(t) \right\} /  \, \{ n^{1/2}\int \rho(t,t)^{1/2} dt \} $)
are fairly low for $n$ small and high for $n$ large.
For instance, when $n=10$, the SNR is 3.16 with the Ornstein-Uhlenbeck process and 6.32 with the Wiener process.
A SNR between 4 and 5 corresponds to a moderate noise level in \cite{hw86}. 
In the  estimation of derivatives, the influence of measurement errors is much stronger.

The simulations were conducted in the R software environment.
All four combinations of $m_1,m_2$ and $\rho_1,\rho_2$ were studied with experimental units $n$
and sampling design size  $N$ varying in $\{10,50,100 \}$.
Different estimation targets ($m_i^{(\nu)}, \, \nu=0,1$) and local polynomial estimators were considered
($p=0,1$ for $\nu=0$, i.e. local constant and linear fits, and $p=1,2$ for $\nu=1$, i.e. local linear and quadratic fits).
In each case, 1000 instances of model \eqref{model} were simulated.
The kernel $K$ was a truncated Gaussian density and the global bandwidths $h$ considered were the exact optimal bandwidth $h_{ex}$ (obtained by minimizing the true IMSE with $w\equiv 1$; see Remark \ref{optimization IMSE}), the asymptotic optimal bandwidth $h_{as}$ (see Corollary \ref{local opt h} and Remark \ref{optimization IMSE}), and the cross-validated bandwidth $h_{cv}$.

Some of the extensive simulation results are presented in Tables \ref{table1}--\ref{table5}. In each table, columns 3-4 are the exact and asymptotic optimal bandwidths; column 5 is the median cross-validated bandwidth over 1000 simulations;
columns 6-7-8 are the median $L^2$ estimation errors $\int_{0}^1 ( \hat{m}_\nu(x)-m^{(\nu)}(x))^2 dx$ 
(first and third quantiles are between brackets) with the exact optimal, asymptotic optimal, and cross-validated bandwidth over the 1000 simulations.
As the estimation errors are strongly right-skewed and feature outliers, 
these errors are described in terms of quantiles rather than mean and standard deviation.

We first comment the results on the estimation of $m$.
Looking at Tables \ref{table1} and \ref{table2} (local linear estimation of $m_1$ with covariance $\rho_1$ or $\rho_2$),
it appears that the bandwidths $h_{ex}, \ h_{as}, \ h_{cv}$
are very close and yield similar performances for almost all $n,N$.
However, note in Table \ref{table1} that $h_{as}$ yields smaller performances when $n=50,100$ and $N=10$,
which can be expected since this  bandwidth is only optimal for large $N$.
Also, local constant estimation of $m_1$ or $m_2$ (not displayed here) yields very similar results to local linear estimation.
In Table \ref{table3} (local linear estimation of $m_2$ with covariance $\rho_2$), the bandwidths $h_{ex}$ and $h_{cv}$ can be infinite for $n=10$.
This remarkable fact has two reasons.
First, the shape of $m_2$ is close to linear, which means that when there are few design points,
increasing the bandwidth $h$ of the estimator only increases its bias marginally (the local linear estimator is unbiased for linear functions). 
On the other hand, under the Ornstein-Uhlenbeck noise,
the estimator's variance reduces drastically as $h$ increases
(much more so than under the Wiener noise).
This can be seen in Theorem 1 where the corrective term in the variance is $-\frac{\alpha(x)}{2}h \iint |u-v| K(u)K(v)dudv $,
with $\alpha(x) =\rho^{(0,1)}(x,x^-)-\rho^{(0,1)}(x,x^+)=30$ for the Ornstein-Uhlenbeck covariance $\rho_2$ and only $\alpha(x)=1$ for the Wiener covariance $\rho_1$.

\bigskip

\begin{center}
$>>>>$ \qquad \quad INSERT TABLES 1-2-3 ABOUT HERE  \quad \qquad $<<<<$
\end{center}

\bigskip

\begin{table}[d]
\begin{center}
\begin{tabular}{ c c | c c c  |  c c c l }

$n$ & $N$ & $h_{ex}$ & $h_{as}$  & $h_{cv}$   & $L^2_{ex}$  & $L^2_{as}$  & $L^2_{cv}$ \\
\hline
10&  10  &  0.07 &   0.06 &  0.08  &      0.031     (0.015-0.061)&0.031     (0.015-0.062)&0.032     (0.016-0.062) \\
  10  &50   & 0.07  &  0.06 &   0.07 &     0.026     (0.012-0.063)    &  0.026     (0.012-0.063)   &   0.027     (0.012-0.065) \\
  10 &100  &  0.07 &   0.06 &   0.07 &      0.028     (0.012-0.060)   &   0.029     (0.012-0.061)  &    0.029     (0.012-0.062) \\
  50  &10   & 0.06   & 0.04 &  0.06 &      0.010     (0.006-0.017)    &  0.018     (0.014-0.026)   &   0.010     (0.006-0.017) \\
  50 & 50   & 0.03  &  0.04 &  0.03   &   0.006     (0.002-0.013)    &  0.006     (0.002-0.013)    &  0.006     (0.002-0.013) \\
  50 &100   & 0.03  &  0.04  & 0.03   &   0.006     (0.002-0.012)    &  0.006     (0.002-0.012)   &   0.006     (0.002-0.012) \\
 100 & 10   & 0.06  &  0.03 &  0.06    &  0.008     (0.005-0.011)    &  0.016     (0.013-0.020)    &  0.008     (0.005-0.011) \\
 100  &50   & 0.03  &  0.03  & 0.03    &  0.003     (0.001-0.006)    &  0.003     (0.001-0.006)    &  0.003     (0.001-0.006) \\
 100 &100   & 0.03  &  0.03  & 0.03  &    0.003     (0.001-0.006)    &  0.003     (0.001-0.006)   &   0.003     (0.001-0.006)
\end{tabular}
\end{center}
\vspace*{2mm}
\caption{Local linear estimation of  $m_1$ with Wiener process noise. }
\label{table1}
\end{table}

\begin{table}[d]
\begin{center}
\begin{tabular}{ c c | c c c  |  c c c l }
$n$ & $N$ & $h_{ex}$ & $h_{as}$  & $h_{cv}$   & $L^2_{ex}$  & $L^2_{as}$  & $L^2_{cv}$ \\
\hline
  10 & 10&0.17&0.15 &0.18 & 0.050 (0.031-0.078)&  0.050 (0.031-0.079)&  0.051 (0.032-0.078) \\
  10 & 50&0.14&0.15  & 0.15&  0.035 (0.021-0.057)&  0.035 (0.021-0.056)&  0.043 (0.025-0.068) \\
  10 &100&0.14&0.15  & 0.16&  0.033 (0.021-0.054)&  0.033 (0.021-0.053)&  0.041 (0.026-0.065) \\
  50 & 10&0.09&0.09  & 0.08&  0.016 (0.010-0.025)&  0.016 (0.010-0.025)&  0.016 (0.011-0.025) \\
  50 & 50&0.08&0.09  & 0.08&  0.009 (0.006-0.014)&  0.009 (0.006-0.014)&  0.010 (0.007-0.015) \\
  50 &100&0.08&0.09  & 0.08&  0.009 (0.006-0.014)&  0.009 (0.006-0.014)&  0.010 (0.007-0.015) \\
 100  &10&0.07&0.07  & 0.08&  0.011 (0.007-0.016)&  0.011 (0.007-0.016)&  0.011 (0.007-0.016) \\
 100 & 50&0.06&0.07  & 0.06&  0.005 (0.003-0.007)&  0.005 (0.003-0.007)&  0.005 (0.004-0.008) \\
 100 &100&0.06&0.07  & 0.06&  0.005 (0.003-0.007)&  0.005 (0.003-0.007)&  0.005 (0.004-0.008)
\end{tabular}
\end{center}
\vspace*{2mm}
\caption{Local linear estimation of  $m_1$ with Ornstein-Uhlenbeck process noise. }
\label{table2}
\end{table}

\begin{table}[d]
\begin{center}
\begin{tabular}{ c c | c c c  |  c c c l }
$n$ & $N$ & $h_{ex}$ & $h_{as}$  & $h_{cv}$   & $L^2_{ex}$  & $L^2_{as}$  & $L^2_{cv}$ \\
\hline
10 & 10   & $\infty$   & 0.28  & $\infty$  &    0.026     (0.017-0.044)&      0.034     (0.020-0.059)&      0.029     (0.018-0.050) \\
10 & 50   & $\infty$  &  0.28  & $\infty$    &  0.024     (0.015-0.039)&      0.025     (0.016-0.040)&      0.027     (0.016-0.044) \\
10& 100 &   $\infty$  &  0.28&   $\infty$   &   0.025     (0.015-0.041)&      0.026     (0.016-0.043)&      0.029     (0.017-0.048) \\
  50  &10   & 0.16  &  0.16  & 0.20   &   0.010     (0.006-0.016)&      0.010     (0.006-0.016)&      0.012     (0.007-0.242) \\
  50 & 50  &  0.12  &  0.16  & 0.14     & 0.008     (0.005-0.012)&      0.008     (0.005-0.012)&      0.010     (0.006-0.018) \\
  50 &100  &  0.12  &  0.16 &  0.14   &   0.008     (0.005-0.012)&      0.008     (0.005-0.012)&      0.010     (0.006-0.018) \\
 100 & 10  &  0.12   & 0.13  & 0.14   &   0.006     (0.004-0.009)&      0.006     (0.004-0.009)&      0.006     (0.004-0.009) \\
 100  &50  &  0.09  &  0.13 &  0.10   &   0.004     (0.003-0.006)&      0.005     (0.003-0.007)&      0.005     (0.003-0.007) \\
 100 &100  &  0.09  &  0.13  & 0.10  &    0.004     (0.003-0.006)&      0.004     (0.003-0.007)&      0.005     (0.003-0.007)
\end{tabular}
\end{center}
\vspace*{2mm}
\caption{Local linear estimation of  $m_2$ with Ornstein-Uhlenbeck process noise. }
\label{table3}
\end{table}

\begin{figure}[d]
\begin{center}
\includegraphics[scale=.75]{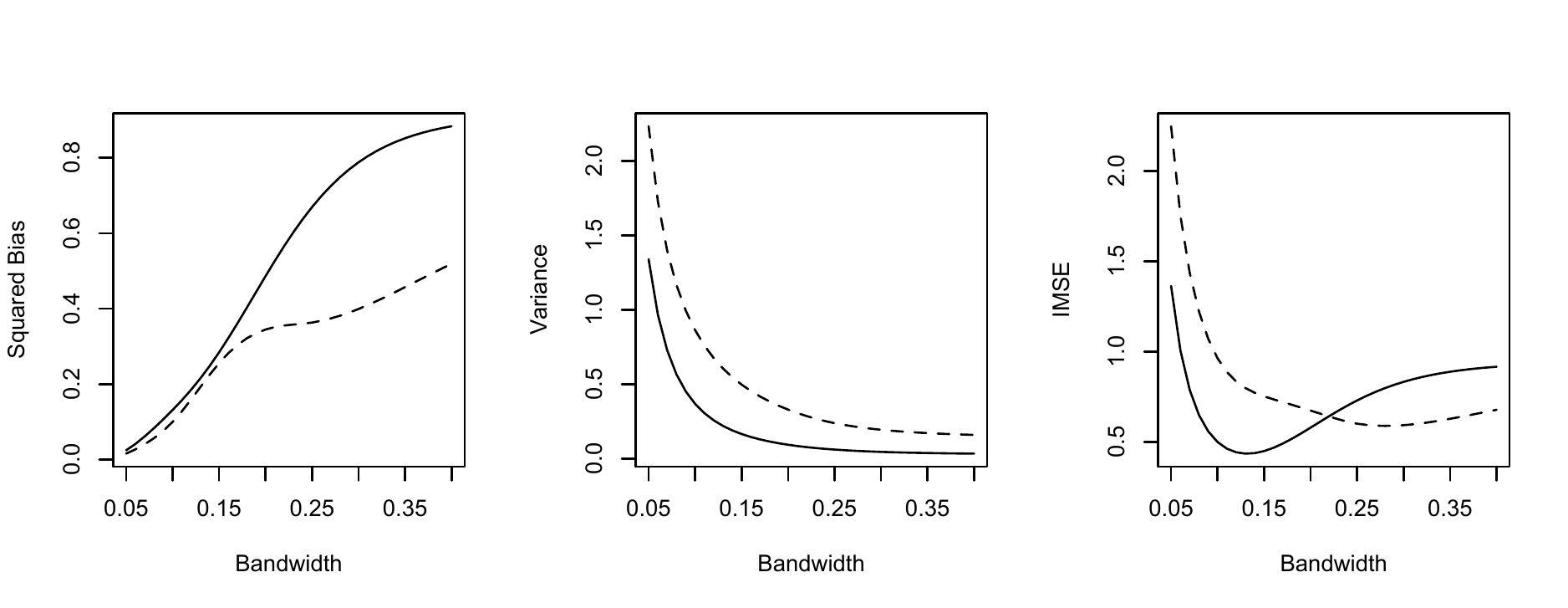}
\caption{Comparison of local linear (solid line) and local quadratic fitting (dashed line) for the estimation of derivatives.
The estimation target is $m_2'$ and the covariance function is $\rho_2$, with $n=N=50$ in \eqref{model}.
  }
\label{linvsquad}
\end{center}
\end{figure}

We now turn to the results concerning the estimation of the derivative $m'$ in Tables 4 and 5 along with Figure 2.
Over the simulations, the use of $h_{ex}$ appears to sensibly reduce the $L^2$-error in comparison to $h_{as}$ and $h_{cv}$.
For the local linear estimation of $m'$, the reduction is $13\%$ and $16\%$, respectively
(median reduction in $L^2$-error over all combinations of $n,N,m_i,$ and $\rho_i$).
The higher performance of $h_{ex}$ over $h_{as}$ and $h_{cv}$ (and of $h_{as}$ over $h_{cv}$ when $N\ge 50$)
can be observed in Tables \ref{table4}-\ref{table5} in the case of the regression $m_1$ and covariance $\rho_1$.
In fact, similar comparisons hold for all choices of $m$ and $\rho$.
For the local quadratic estimation of $m'$, the use of $h_{ex}$ and $h_{cv}$ reduce the $L^2$-error by respectively $45\%$ and $30\%$, in comparison to $h_{as}$ (over all combinations of $n,N,m_i,$ and $\rho_i$). 
It is noteworthy that $h_{as}$ is systematically smaller than $h_{ex}$
(the difference between the two bandwidths is larger when estimating $m'$ than when estimating $m$) 
but the two bandwidths $h_{as}$ and $h_{ex}$ are closer for $p=1$ (local linear fit) than for $p=2$ (local quadratic). \\
Comparing local linear to local quadratic estimation, the latter can considerably reduce the bias at the expense of increasing the variance,
which is a consequence of adding an extra (quadratic) parameter in the local fit.
Which order of local polynomial fit achieves better performances in a given scenario depends on the balance between bias and variance.
It can be seen from Tables \ref{table4} and \ref{table5} that when $h_{ex}$ is used in the simulations, the local quadratic estimator yields sensibly better results than the local linear when the target is $m_1'$
(due to the fairly high curvature of $m_1'$ which makes the bias large in comparison to the variance).
The situation is however reversed with the target $m_2'$ (that has relatively low curvature), as shown in Figure \ref{linvsquad}.
In this figure,  the local quadratic estimator has a slightly smaller (squared and integrated) bias than the local linear for small $h$
(see left panel). On the other hand, the local linear estimator has much smaller variance than the quadratic for all $h$ (middle panel).
Overall, the optimal IMSE is smaller for the local linear estimator and the optimal bandwidths are quite different (right panel), $h_{opt}\approx 0.13$ for the linear one and $h_{opt}\approx 0.30$ for the quadratic one.

\bigskip

\begin{center}
$>>>>>>$ \qquad \quad INSERT TABLES 4-5 ABOUT HERE  \quad \qquad $<<<<<<$
\end{center}

\bigskip

\begin{center}
$>>>>>>$ \qquad \quad INSERT FIGURE 2 ABOUT HERE  \quad \qquad $<<<<<<$
\end{center}

\bigskip

\begin{table}[d]
\begin{center}
\begin{tabular}{ c c | c c c  |  c c c l }
$n$ & $N$ & $h_{ex}$ & $h_{as}$  & $h_{cv}$   & $L^2_{ex}$  & $L^2_{as}$  & $L^2_{cv}$ \\
\hline
10&10&0.07&0.04&0.08& 1.97 (1.52-2.61)&3.97 (3.72-4.30)&2.05 (1.59-2.72)\\
10&50&0.05&0.04&0.07&0.84 (0.62-1.13)&0.88 (0.64-1.15)&1.07 (0.77-1.50)\\
10&100&0.05&0.04&0.07&0.84 (0.60-1.13)&0.90 (0.64-1.15)&1.13 (0.79-1.60)\\
50&10&0.06&0.03&0.06&1.63 (1.42-1.86)&4.18 (4.14-4.24)&1.65 (1.43-1.87)\\
50&50&0.03&0.03&0.03&0.29 (0.22-0.37)&0.28 (0.22-0.37)&0.30 (0.23-0.39)\\
50&100&0.03&0.03&0.03&0.26 (0.20-0.33)&0.26 (0.20-0.33)&0.28 (0.21-0.36)\\
100&10&0.06&0.03&0.06&1.60 (1.45-1.77)&3.73 (3.69-3.76)&1.61 (1.47-1.78)\\
100&50&0.02&0.02&0.03&0.18 (0.14-0.23)&0.18 (0.14-0.22)&0.18 (0.14-0.22)\\
100&100&0.02&0.02&0.03&0.17 (0.13-0.20)&0.17 (0.13-0.20)&0.17 (0.13-0.21)
\end{tabular}
\end{center}
\vspace*{2mm}
\caption{ Local linear estimation of  $m_1'$ with Wiener process noise. }
\label{table4}
\end{table}

\begin{table}[d]
\begin{center}
\begin{tabular}{ c c | c c c  |  c c c l }
$n$ & $N$ & $h_{ex}$ & $h_{as}$  & $h_{cv}$   & $L^2_{ex}$  & $L^2_{as}$  & $L^2_{cv}$ \\
\hline
 10 & 10   &  0.11   &  0.04   & 0.08    &   0.96      (0.59-1.71)& 23.0     (11.1-50.1)&       3.41      (1.39-6.89) \\
 10 & 50    & 0.09   &  0.04  &  0.07     &  0.52      (0.33-0.76)&       0.92      (0.67-1.21)&       0.64      (0.42-0.94) \\
  10 &100  &   0.09 &    0.04  &  0.07   &    0.52      (0.33-0.78)&       0.89      (0.64-1.17)&       0.63      (0.42-0.91) \\
  50 & 10   &  0.09  &   0.03  &  0.06    &   0.52      (0.30-0.81)&    1940    (536-4894)&       8.27      (7.76-8.90) \\
  50 & 50   &  0.06   &  0.03  & 0.03     &  0.15      (0.10-0.21)&       0.26      (0.19-0.35)&       0.23      (0.17-0.32) \\
  50 &100  &   0.06  &   0.03  &  0.03   &    0.14      (0.09-0.19)&       0.23      (0.18-0.29)&       0.20      (0.15-0.26) \\
 100  &10   &  0.09   &  0.03  &  0.06    &   0.48      (0.32-0.69)&      23.8     (15.7-33.0)&       8.20      (7.82-8.64) \\
 100 & 50   &  0.05  &   0.02  &  0.03     &  0.08      (0.06-0.12)&       0.15      (0.12-0.20)&       0.15      (0.12-0.20) \\
 100 &100   &  0.05   &  0.02  &  0.03    &   0.08      (0.06-0.11)&       0.13      (0.10-0.17)&       0.13      (0.10-0.16)
\end{tabular}
\end{center}
\vspace*{2mm}
\caption{ Local quadratic estimation of  $m_1'$ with Wiener process noise. }
\label{table5}
\end{table}

\section{Discussion}

We have examined in this paper the local polynomial estimation of a regression function and its derivatives in the context of functional data.
Our main theoretical contribution has been to derive second-order asymptotic expansions for the bias and variance of the estimator based on a sampling density not necessarily uniform.
These expansions give qualitative insights on the large-sample behavior of the estimators and highlight in particular
how the covariance and the choice of the bandwidth affect the estimator's variance.
Our result fills an important gap in the literature as, to this date,
no asymptotic theory seems available on the estimation of regression derivatives with functional data under correlated errors.
This topic is relevant in practice since for many functional data sets,
essential information may be carried by derivatives of the observed curves.
Note that our results may be extended to the multivariate regression setup and also to noisy functional data.

We have applied our main result to determine optimal sampling densities and bandwidths 
that can be estimated in practice by plugin methods.
To examine the potential benefits of a plugin method for bandwidth selection,
we have compared numerically the performances of local polynomial estimators based on the asymptotic optimal bandwidth $h_{as}$ (nessary for the plugin method),
the exact optimal bandwidth $h_{ex}$, and the cross-validated bandwidth $h_{cv}$ of \citet{rs91}. 
The simulations indicate that a plugin method could be an interesting alternative to cross-validation
for data sets with moderate to large numbers of observation points (which is typically the case for functional data),
especially for estimating the derivatives of the regression function $m$.
Developing a plugin method would however require to estimate
the partial derivatives of the covariance function $\rho$ and some higher-order derivative of $m$.
Another outcome of the simulations is that although cross-validation is not meant for derivative estimation,
estimators based on $h_{cv}$ generally give satisfactory results both when the target is  $m$ and $m'$.
Finally, our simulations suggest the use of local linear fits both for estimating $m$ (rather than local constant) and $m'$ (rather than local quadratic),
as these estimators are more stable (especially for small sample sizes $N$) and give reasonable estimates in most situations.

Finally, we have established the asymptotic normality of the local polynomial estimator in the pointwise sense.
This result can be applied in various inference procedures.
By following the arguments of \cite{d10}, 
a stronger asymptotic normality result can be obtained in the space of continuous functions equipped with the sup-norm. This allows  to conduct simultaneous inference on the regression derivatives.


\section*{Appendix: Proofs}

Throughout the proofs, the dependence of vectors and matrices on $N$
is denoted explicitly to clarify the arguments. Also, to fix ideas, the compact support of $K$
is taken to be $[-1,1]$ without loss of generality.

\subsection*{Proof of Theorem 1: bias term}

Let us write $\mathbf{m}_N= (m(x_1),\ldots, m(x_N))'$ and define the
$(p+1)\times(p+1)$
matrix $\mathbf{S}_N = N^{-1} \mathbf{X}_N' \mathbf{W}_N \mathbf{X}_N $ with $(k,l)$th element
($0\le k,l \le p$)
given by
\begin{equation*}
s_{k+l,N}=  \frac{1}{Nh} \sum_{j=1}^N (x_j-x)^{k+l} K\left(\frac{x_j-x}{h}\right) . 
\end{equation*}

It follows from \eqref{def estimateur}
that
\begin{equation}
\mathbb{E} (\boldsymbol{\widehat{\beta}}_N (x))=  N^{-1} \mathbf{S}_N^{-1}\mathbf{ X}_N' \mathbf{W}_N \mathbf{m}_N .
\end{equation}

With (A3), a Taylor expansion of $m(x_j)$ at the order $(p+2)$ yields
\begin{equation*}
m(x_j)= m(x) + (x_j - x)m'(x) +\ldots + \frac{ (x_j-x)^{p+2}}{(p+2)!} m^{(p+2)}(x)+ o\left( (x_j-x)^{p+2}\right)
\end{equation*}
and thus
$$
\mathbf{m}_N= \mathbf{X}_N  \boldsymbol{\beta}(x) + \beta_{p+1} \left(\begin{array}{c}
(x_1-x)^{p+1} \\
\vdots \\
(x_N- x)^{p+1}
\end{array}
\right)
+
( \beta_{p+2} +o(1) ) \left(\begin{array}{c}
(x_1-x)^{p+2} \\
\vdots \\
(x_N- x)^{p+2}
\end{array}
\right).
$$
Hence the bias in the estimation of $\boldsymbol{\beta}(x)$ is
\begin{equation}\label{bias m_j hat}
\mathbb{E} (\boldsymbol{\widehat{\beta}}_N (x) )- \boldsymbol{\beta}(x) = \beta_{p+1} \mathbf{S}_N^{-1}  \mathbf{c}_N+
\left( \beta_{p+2} +o(1) \right)  \mathbf{S}_N^{-1}
\mathbf{\tilde{c}}_N ,
\end{equation}
where $\mathbf{c}_N= (s_{p+1,N},\ldots, s_{2p+1,N})'$ and $ \mathbf{\tilde{c}}_N=(s_{p+2,N},\ldots, s_{2p+2,N})'$.

With the regularity of the sampling design \eqref{sampling design},
the Lipschitz-continuity of $K$ and the compacity of its support,
some straightforward calculations allow to
approximate the elements of the matrix $\mathbf{S}_N$ by
\begin{equation}\label{approx sn}
s_{k+l,N}= h^{k+l}\, \Big(  \int_{-\infty}^\infty u^{k+l} K(u)f(x+hu)du + \mathcal{O}  ( (Nh)^{-1}  ) \Big) ,
\end{equation}
the $ \mathcal{O}  ( (Nh)^{-1}) $ being the error in the integral approximation of a Riemann sum.

From assumption (A4),  a Taylor expansion of $f(x+hu)$ at order 1 yields
\begin{equation}
s_{k+l,N}=  h^{k+l} \big(\mu_{k+l} f(x) + h \mu_{k+l+1}f'(x) +o(h) \big)
\end{equation}
under the condition $Nh^2\to\infty$, which
in matrix form stands as
\begin{equation}
\mathbf{S}_N=  \mathbf{H}\left(f(x)\mathbf{S} + hf'(x)\mathbf{\tilde{S}} +o(h)\right)\mathbf{H} \,,
\end{equation}
where $\mathbf{H}=diag(1,h,\cdots,h^p)$.
In particular, it holds that
\begin{equation}\label{asymptotic mu_n}
\left\{
\begin{array}{ll}
\mathbf{c}_N& =  h^{p+1} \mathbf{H} \big(f(x)\mathbf{c} + \left( hf'(x)+ o(h) \right) \mathbf{\tilde{c}} \big) , \\
 \mathbf{\tilde{c}}_N &= h^{p+2} \mathbf{H}  \left(f(x)+o(1)\right) \mathbf{\tilde{c}} .
 \end{array} \right.
\end{equation}

Now, due to the fact that $(\mathbf{U}+h\mathbf{V})^{-1}= \mathbf{U}^{-1}- h\mathbf{U^{-1}VU^{-1}} +o(h)$
for any two invertible matrices $\mathbf{U,V}$ of compatible dimensions,  we have
\begin{equation}\label{invS}
\mathbf{S}_N^{-1}   =  \mathbf{H}^{-1}\left(\frac{1}{f(x)}\mathbf{S}^{-1}  - h\frac{f'(x)}{f^2(x)}\mathbf{S^{-1}\tilde{S}S^{-1} } + o(h) \right)\mathbf{H}^{-1} .
\end{equation}

Plugging \eqref{asymptotic mu_n}-\eqref{invS} in  \eqref{bias m_j hat} and truncating the expansion to the second order,
the bias expression of Theorem \ref{t1} follows.


\subsection*{Proof of Theorem 1: variance term}

Define the $N\times N$ matrix $\boldsymbol{\Sigma }_N =  \left(\rho(x_i,x_j) \right)$ and the $(p+1)\times(p+1)$ matrix $\mathbf{S}_N^*= N^{-2} \mathbf{  X}_N' \mathbf{W}_N
\boldsymbol{\Sigma}_N \mathbf{W}_N \mathbf{ X}_N$.
 Noting that $\mathrm{Var}(\bar{\mathbf{Y}})= n^{-1}\boldsymbol{\Sigma}_N $
and considering \eqref{def estimateur},
it can be seen that
\begin{equation}\label{variance matrix form}
\mathrm{Var} (\boldsymbol{\widehat{\beta}}_N (x))
= n^{-1} \mathbf{S}_N^{-1}\mathbf{S}_N^* \mathbf{S}_N^{-1}\, .
\end{equation}

The asymptotic behavior of the matrix $\mathbf{ S}_N^\ast$ is given by the following lemma.
\begin{lem}\label{Sn ast}
Assume (A1)-(A5) for a given $x\in(0,1)$.  Then as $n,N \to \infty $,
\begin{align*}
\mathbf{S}_N^* & =   \mathbf{H} \left\{ \phi(x,x) \mathbf{S}^{\ast} + h ( \phi ^{(0,1)}(x,x^+)- \phi ^{(0,1)}(x,x^-)) \mathbf{A} \right. \\
& \qquad \qquad \quad \left. +h( \phi ^{(0,1)}(x,x^+)+ \phi ^{(0,1)}(x,x^-))\mathbf{B} +o(h) \right\} \mathbf{H}
\end{align*}
with $\mathbf{A,S^\ast}$ being defined in Section 3, $\mathbf{B}=\left(\frac{1}{2}(\mu_{k+1}\mu_l + \mu_{k}\mu_{l+1})\right)$, and   $\phi(y,z) = \rho(y,z) f(y) f(z)$.
\end{lem}

Plugging Lemma \ref{Sn ast} and \eqref{invS} in \eqref{variance matrix form}, we have
\begin{equation}
\begin{split}
 n f(x)^2 \mathbf{H} \mathrm{Var} (\boldsymbol{\widehat{\beta}}(x)) \mathbf{H}
&  = \phi(x,x)  \, \mathbf{S}^{-1} \mathbf{S^\ast}\mathbf{S}^{-1}   +o(h)  \\
&- h \phi(x,x)\, \frac{f'(x)}{f(x)} \,
 ( \mathbf{S}^{-1}\tilde{\mathbf{S}}\mathbf{S}^{-1}\mathbf{S}^{\ast} \mathbf{S}^{-1}
+ \mathbf{S}^{-1}\mathbf{S}^{\ast}\mathbf{S}^{-1} \tilde{\mathbf{S}} \mathbf{S}^{-1}) \\
&\quad + h \, (\phi ^{(0,1)}(x,x^+) - \phi ^{(0,1)}(x,x^-) )\, \mathbf{S}^{-1} \mathbf {AS}^{-1} \\
& \quad + h \, (\phi ^{(0,1)}(x,x^+) + \phi ^{(0,1)}(x,x^-) )\, \mathbf{S}^{-1} \mathbf {BS}^{-1}.
\end{split}
\end{equation}
Note that the $o(h)$ above stands for a matrix whose coefficients are negligible compared to  $h$ as $h\to 0$.

Expressing $\phi ^{(0,1)}(x,x^\pm)$ in terms of $\rho^{(0,1)}(x,x^\pm)$, we get
\begin{equation}
\phi^{(0,1)}(x,x^\pm)=f(x)f'(x)\rho(x,x) + f^2(x)\rho^{(0,1)}(x,x^\pm)
\end{equation}
and then
\begin{equation}\label{var0}
\begin{split}
 n \mathbf{H} & \mathrm{Var} (\boldsymbol{\widehat{\beta}}(x))  \mathbf{H}
= \rho(x,x) \mathbf{S}^{-1}\mathbf{S^\ast S}^{-1}  + o(h) \\
&  - h \rho(x,x)\, \frac{f'(x)}{f(x)} \,
 ( \mathbf{S}^{-1}\tilde{\mathbf{S}}\mathbf{S}^{-1}\mathbf{S}^{\ast} \mathbf{S}^{-1}
+ \mathbf{S}^{-1}\mathbf{S}^{\ast}\mathbf{S}^{-1} \tilde{\mathbf{S}} \mathbf{S}^{-1}) \\
& + h \, (\rho^{(0,1)}(x,x^+) - \rho^{(0,1)}(x,x^-) )\, \mathbf{S}^{-1} \mathbf {AS}^{-1} \\
&  + h\left(2 \frac{f'(x)}{f(x)}\rho(x,x) +  (\rho ^{(0,1)}(x,x^+) + \rho ^{(0,1)}(x,x^-) )\right) \mathbf{S}^{-1} \mathbf {BS}^{-1} .
\end{split}
\end{equation}

The variance expression can further be simplified due to the fact that
\begin{equation}\label{var1}
\left\{ \begin{array}{l}
 \mathbf{e}_{\nu}' \mathbf{S}^{-1}\tilde{\mathbf{S}}\mathbf{S}^{-1}\mathbf{S}^{\ast} \mathbf{S}^{-1}  \mathbf{e}_{\nu} = 0 \\
 \mathbf{e}_{\nu}' \mathbf{S}^{-1}  \mathbf{B}  \mathbf{S}^{-1} \mathbf{e}_{\nu} = 0
 \end{array} \right.
\end{equation}
for all $\nu=0,\ldots,p$. To see this, we need to examine in detail  the above matrices.
By the symmetry of $K$, $\mathbf{S}=(\mu_{k+l})$ has its $(k,l)$th entry equal to zero if $k,l$ are of different parity.
The same property can be established for $\mathbf{S}^{-1}$ by standard cofactor arguments.
For $\tilde{\mathbf{S}}=(\mu_{k+l+1})$, the $(k,l)$th entry is zero if $k,l$ are of the same parity.
For $\mathbf{S}^\ast = (\mu_k \mu_l)$, the sparsity is even stronger: all the rows, columns, and subdiagonals of odd order
(recall that the indexing starts at 0) have their entries equal to zero. With some matrix algebra, one can check that the matrices
 $\mathbf{S}^{-1}\tilde{\mathbf{S}}\mathbf{S}^{-1}$ and $\mathbf{S}^{\ast} \mathbf{S}^{-1}$
 have the same sparsity structures as $\tilde{\mathbf{S}}$ and $\mathbf{S}^\ast$, respectively.
 It is then easy to obtain the first part of  \eqref{var1}. The second part is derived along the same lines.
It suffices to notice, on the one hand, that
 $ \mathbf{e}_{\nu}' \mathbf{S}^{-1}   \mathbf{BS}^{-1} \mathbf{e}_{\nu}$  can be written as the double sum $ \sum_{k,l} \left[\mathbf{S}^{-1} \mathbf{e}_{\nu}\right]_k  \left[\mathbf{S}^{-1} \mathbf{e}_{\nu}\right]_l   B_{kl}$ over the indexes $k,l$ having the same parity as $\nu$. On the other hand, $B_{kl} = \mu_k \mu_{l+1}+\mu_{k+1}\mu_l =0$ for $k,l$ both even or both odd ($\mu_k=\mu_l=0$ if $k,l$ odd and $\mu_{k+1}=\mu_{l+1}=0$ if $k,l$ even). Combining these two facts yields the asymptotic result for the variance term.

Finally, we deduce from \eqref{var0} and \eqref{var1} that
\begin{align*}
n \mathrm{Var} (\hat{m}_{\nu}(x)) & = n (\nu!)^2\, \mathbf{e}_\nu' \, \mathrm{Var} (\boldsymbol{\widehat{\beta}}(x))\,  \mathbf{e}_\nu \\
& =  (\nu!)^2 h^{-2\nu}\, \rho(x,x)\, \mathbf{e}_\nu'  \mathbf{S}^{-1}\mathbf{S^\ast S}^{-1}  \mathbf{e}_\nu + o(h^{-2\nu +1 })\\
&\quad  + (\nu!)^2 h^{-2\nu +1 }\, (\rho^{(0,1)}(x,x^+) - \rho^{(0,1)}(x,x^-)) \, \mathbf{e}_\nu' \mathbf{S}^{-1} \mathbf{A}\mathbf{S}^{-1} \mathbf{e}_\nu ,
\end{align*}
which completes the proof of  Theorem 1. $\square$


\subsection*{ Proof of Lemma  \ref{Sn ast}.}

By the same arguments used to approximate the matrix $\mathbf{S}_N$ with integrals in \eqref{approx sn}, one can use the regularity (A5) of the covariance function $\rho$ to show that the elements  of $ \mathbf{S}^\ast_{N}$ satisfy 
 \begin{align}\label{s-ast kl}
& s^\ast_{kl,N}= (Nh)^{-2} \sum_{i=1}^N \sum_{j=1}^N (x_i-x)^k (x_{j}-x)^{l}K\left(\frac{x_i-x}{h}\right)K\left(\frac{x_j-x}{h}\right)\rho(x_i,x_{j}) \nonumber \\
&=    h^{-2}
\iint_{[-1,1]^2}(u-x)^k (v-x)^{l}K\left(\frac{u-x}{h}\right)K\left(\frac{v-x}{h}\right) \rho(u, v)f(u)f(v) dudv  \nonumber \\
& \qquad + \mathcal{O}\left(\frac{ h^{k+l}}{Nh}\right)  \nonumber \\
&=   h^{k+l} \iint_{[-1,1]^2} u^k v^{l} \phi(x+hu, x+hv) K(u)K(v)dudv +o( h^{k+l+1}),
\end{align}
assuming that $(Nh)^{-1}=o(h)$, i.e. $Nh^2 \to \infty$.

Using Taylor expansions together with (A4)-(A5), one can show that
\begin{equation*}
\phi(x+hu,x+hv) = \phi(x,x) + hu \phi^{(0,1)}(x,x^{-}) + hv \phi^{(0,1)}(x,x^+) +o(h)
\end{equation*}
 for all $0<u<v<1$.
This expansion is obtained by introducing a pivotal point $(x+hu,x)$ or $(x,x+hv)$ such that the lines
 connecting this point to $(x+hu,x+hv)$ and $(x,x)$ do not cross the main diagonal of $[0,1]^2$.
 One can then safely perform Taylor expansions along the connecting lines,
 knowing that $\phi$ is differentiable on each side of the diagonal.
 The above expansion also relies on the identities
 $\phi^{(1,0)}(x^{+},x) =\phi^{(0,1)}(x,x^{+}) = \phi^{(0,1)}(x^-,x)$
 and $\phi^{(1,0)}(x^{-},x) = \phi^{(0,1)}(x,x^-)= \phi^{(0,1)}(x^+,x)$
 (thanks to the symmetry of $\phi$ and the continuity of the first
partial derivatives of $\phi$ on either side of the diagonal).
 By symmetry considerations, it then holds for all $u,v\in[0,1]$ that
\begin{equation}\label{taylor expansion 1}
\begin{split}
 \phi(x+hu,x+hv) & = \phi(x,x)  + h\, (u\wedge v) \, \phi^{(0,1)}(x,x^{-}) \\
 & \qquad + h \, (u \vee v) \,\phi^{(0,1)}(x,x^+) +o(h).
\end{split}
\end{equation}

Using the fact that $(u\wedge v)+(u\vee v)=u+v$ and $(u\vee v)- (u\wedge v)= |u-v|$,
writing $\phi^{(0,1)}(x,x^\pm) = \frac{1}{2}\,(\phi^{(0,1)}(x,x^+)+\phi^{(0,1)}(x,x^-)) \pm \frac{1}{2} \,( \phi^{(0,1)}(x,x^+)-\phi^{(0,1)}(x,x^-))$, 
one concludes, by the dominated convergence theorem, that
\begin{align*}
s_{kl,N}^\ast  &= h^{k+l} \iint_{[-1,1]^2} u^k v^l K(u) K(v) \phi(x+hu,x+hv)dudv + o (h^{k+l+1}) \\
& = h^{k+l} \Big\{ \phi(x,x) \mu_k \mu_l + \frac{h}{2} \left( \phi^{(0,1)}(x,x^{+})+\phi^{(0,1)}(x,x^{-})   \right)  \left( \mu_{k+1}\mu_l + \mu_k \mu_{l+1}\right) \\
&  \quad + \frac{h}{2} \left( \phi^{(0,1)}(x,x^{+}) - \phi^{(0,1)}(x,x^{-})   \right) \iint_{[-1,1]^2} |u-v| u^k v^l K(u) K(v)du dv\Big\} \\
& \quad + o(h^{k+l+1})  . \quad \square
 \end{align*}


\subsection*{ Proof of Theorem  \ref{t2}.}

This result is obtained along the lines of the proof of Theorem \ref{t1}. More precisely,
it suffices to push the matrix expansions of  $\mathbf{S}_N^{-1}$ in
 \eqref{invS} and $\mathbf{S}_N^\ast$ in Lemma \ref{Sn ast} to a higher order $d$.
First, since $f\equiv 1$, it is easily seen that $\mathbf{S}_N = \{ 1+ o(h^d)\} \mathbf{HSH}$
provided that $Nh^{d+1}\to\infty$. Therefore,  \eqref{invS} simply extends in
$\mathbf{S}_N^{-1} =  \{1+ o(h^d)\} \mathbf{H}^{-1}\mathbf{S}^{-1}\mathbf{H}^{-1}$.
Second, if the covariance $\rho$ is $d$ times differentiable at $(x,x)$,
then a Taylor expansion of order $d$ can be performed for $\rho(x+hu,x+hv)$,
followed by an application of the dominated convergence theorem over $[-1,1]^2$ as $h\to 0$.
For $d=4$, we get for instance (see the proof of Lemma  \ref{Sn ast}):
\begin{align}\label{SN ast 2}
s_{kl,N}^\ast  &= h^{k+l} \iint_{[-1,1]^2} u^k v^l K(u) K(v) \rho(x+hu,x+hv)dudv + o (h^{k+l+4}) \nonumber \\
& = h^{k+l} \Big\{ \rho(x,x) \mu_k \mu_l +  h \rho^{(0,1)}(x,x)  \left( \mu_{k+1}\mu_l + \mu_k \mu_{l+1}\right) \nonumber\\
&  \qquad + h^2 \left( \rho^{(0,2)}(x,x) \frac{\mu_{k+2} \mu_{l} + \mu_{k} \mu_{l+2}}{2!} +  \rho^{(1,1)}(x,x) \mu_{k+1}\mu_{l+1}\right) \nonumber \\
& \qquad  + h^3  \left( \rho^{(0,3)}(x,x) \frac{\mu_{k+3} \mu_{l} + \mu_{k} \mu_{l+3}}{3!} +  \rho^{(1,2)}(x,x) \frac{ \mu_{k+2}\mu_{l+1} +\mu_{k+1}\mu_{l+2}}{2!} \right) \nonumber \\
& \qquad  + h^4  \left( \rho^{(0,4)}(x,x) \frac{\mu_{k+4} \mu_{l} + \mu_{k} \mu_{l+4}}{4!} +  \rho^{(1,3)}(x,x) \frac{ \mu_{k+3}\mu_{l+1} +\mu_{k+1}\mu_{l+3}}{3!}
\right. \nonumber \\
&\qquad \qquad  \qquad \left. +  \rho^{(2,2)}(x,x) \frac{\mu_{k+2} \mu_{l+2}}{2!\,  2!}  \right) + o(h^4) \Big\}  .
 \end{align}

The arguments used in Theorem \ref{t1} relative to the sparsity structure of $\mathbf{S}^{-1}$ and the limit matrix of $\mathbf{S}_N^\ast$ still apply here.
In a nutshell, the matrices of the form $(\mu_{k+a}\mu_{l+b})$ in \eqref{SN ast 2} that do contribute  to the limit variance of $\hat{m}_\nu(x)$
are those for which both $\nu+a$ and $\nu+b$ are even. (This corresponds to the nonzero moments of the kernel $K$.)
Therefore, the terms of order $h$ and $h^3$ inside the brackets of \eqref{SN ast 2}
do not contribute to the limit variance of $\hat{m}_\nu(x)$.
For $\nu$ even, the terms $\mu_{k+1}\mu_{l+1}$ in \eqref{SN ast 2} do not contribute either but
the terms $\mu_k \mu_l$ and $h^2 ( \mu_{k+2} \mu_{l}+ \mu_{k} \mu_{l+2})$ do. An expansion to order $d=2$ is thus sufficient.
 For $\nu$ odd, only the terms $h^2 \mu_{k+1} \mu_{l+1}$ and $h^4 (\mu_{k+3} \mu_{l+1} + \mu_{k+1} \mu_{l+3})$ contribute
to the limit variance of $\hat{m}_\nu(x)$   up to order 4.
 In this case the expansion to order $d=4$ is necessary, as an expansion to order 2 only results in a variance term of order $1/n$ (independent of $h$) when $\nu=1$.
 Theorem \ref{t2} immediately  follows from these arguments. $\square$


\subsection*{ Proof of Lemma  \ref{alpha nonnegative}.}

Starting from the Taylor expansion \eqref{taylor expansion 1} and the subsequent argument in the proof of Lemma \ref{Sn ast},
it can be shown that
\begin{equation}
\begin{split}
\rho(x+hu,x+hv) & = \rho(x,x) +
\frac{h}{2} \left( \rho^{(0,1)}(x,x^{+})+\rho^{(0,1)}(x,x^{-})   \right)  (u+v) \\
& + \frac{h}{2} \left( \rho^{(0,1)}(x,x^{+}) - \rho^{(0,1)}(x,x^{-})   \right) |u-v| + o(h)
\end{split}
\end{equation}
for all $u,v \in [-1,1]^2$ as $h\to 0$.

Let us write $a= \frac{ \rho^{(0,1)}(x,x^{+})+\rho^{(0,1)}(x,x^{-}) }{2}$ and $b=\frac{ \rho^{(0,1)}(x,x^{+})-\rho^{(0,1)}(x,x^{-}) }{2}$
for brevity.
The dominated convergence theorem and (A5)
imply that for any bounded, measurable function $g$ on $[-1,1]$,
\begin{equation}\label{dominated convergence}
\begin{split}
\iint_{[-1,1]^2} & \rho(x+hu,x+hv)g(u)g(v) dudv  \\
& = \rho(x,x) \left( \int_{-1}^{1} g(u)du \right)^2  
 + 2ah    \int_{-1}^1 g(u)du \int_{-1}^1 vg(v)dv \\
& \qquad +bh \iint_{[-1,1]^2} g(u)g(v) |u-v| dudv + o(h) .
\end{split}
\end{equation}

The left handside of \eqref{dominated convergence}
is non-negative since the covariance $\rho$ is a non-negative definite function. By taking $g=\mathrm{Id}_{[-1,1]}$, we have $\int_{-1}^1 g(u)du = 0$ so that
the remaining term $bh \iint_{[-1,1]^2} g(u)g(v) |u-v| dudv $
 in the right handside of \eqref{dominated convergence} is also non-negative. 
Since $\iint_{[-1,1]^2} uv |u-v| dudv = -\frac{8}{15}<0$, this means that $b \le 0$
and hence $\alpha(x)=\rho^{(0,1)}(x,x^-)- \rho^{(0,1)}(x,x^+) \ge 0$.
$\square$


\section*{Acknowledgement}
The authors would like to thank Giles Hooker from Cornell University for his valuable comments that helped enhancing the scope of the paper.


\bibliographystyle{model1b-num-names} 
\bibliography{benh-degras}

\begin{thebibliography}{27}
\expandafter\ifx\csname natexlab\endcsname\relax\def\natexlab#1{#1}\fi
\providecommand{\bibinfo}[2]{#2}
\ifx\xfnm\relax \def\xfnm[#1]{\unskip,\space#1}\fi
\bibitem[{Benhenni and Cambanis(1992)}]{bc92}
\bibinfo{author}{K.~Benhenni}, \bibinfo{author}{S.~Cambanis},
  \bibinfo{title}{Sampling designs for estimating integrals of stochastic
  processes}, \bibinfo{journal}{Ann. Statist.} \bibinfo{volume}{20}
  (\bibinfo{year}{1992}) \bibinfo{pages}{161--194}.
\bibitem[{Benhenni and Rachdi(2006)}]{br06}
\bibinfo{author}{K.~Benhenni}, \bibinfo{author}{M.~Rachdi},
  \bibinfo{title}{Nonparametric estimation of the regression function from
  quantized observations}, \bibinfo{journal}{Comput. Statist. Data Anal.}
  \bibinfo{volume}{50} (\bibinfo{year}{2006}) \bibinfo{pages}{3067--3085}.
\bibitem[{Benhenni and Rachdi(2007)}]{br07}
\bibinfo{author}{K.~Benhenni}, \bibinfo{author}{M.~Rachdi},
  \bibinfo{title}{Nonparametric estimation of average growth curve with general
  nonstationary error process}, \bibinfo{journal}{Comm. Statist. Theory
  Methods} \bibinfo{volume}{36} (\bibinfo{year}{2007})
  \bibinfo{pages}{1173--1186}.
\bibitem[{Cambanis(1985)}]{c85}
\bibinfo{author}{S.~Cambanis}, \bibinfo{title}{Sampling designs for time
  series}, in: \bibinfo{booktitle}{Time series in the time domain},
  volume~\bibinfo{volume}{5} of \textit{\bibinfo{series}{Handbook of
  Statist.}}, \bibinfo{publisher}{North-Holland}, \bibinfo{address}{Amsterdam},
  \bibinfo{year}{1985}, pp. \bibinfo{pages}{337--362}.
\bibitem[{Cardot(2000)}]{car00}
\bibinfo{author}{H.~Cardot}, \bibinfo{title}{Nonparametric estimation of
  smoothed principal components analysis of sampled noisy functions},
  \bibinfo{journal}{J. Nonparametr. Statist.} \bibinfo{volume}{12}
  (\bibinfo{year}{2000}) \bibinfo{pages}{503--538}.
\bibitem[{Degras(2008)}]{d08}
\bibinfo{author}{D.~Degras}, \bibinfo{title}{Asymptotics for the nonparametric
  estimation of the mean function of a random process},
  \bibinfo{journal}{Statist. Probab. Lett.} \bibinfo{volume}{78}
  (\bibinfo{year}{2008}) \bibinfo{pages}{2976--2980}.
\bibitem[{Degras(2010)}]{d10}
\bibinfo{author}{D.~Degras}, \bibinfo{title}{Simultaneous confidence bands for
  nonparametric regression with functional data}, \bibinfo{journal}{Statist.
  Sinica}  (\bibinfo{year}{2010}). \bibinfo{note}{Accepted for publication}.
\bibitem[{Fan and Gijbels(1996)}]{fg96}
\bibinfo{author}{J.~Fan}, \bibinfo{author}{I.~Gijbels}, \bibinfo{title}{Local
  polynomial modelling and its applications}, volume~\bibinfo{volume}{66} of
  \textit{\bibinfo{series}{Monographs on Statistics and Applied Probability}},
  \bibinfo{publisher}{Chapman \& Hall}, \bibinfo{address}{London},
  \bibinfo{year}{1996}.
\bibitem[{Fan et~al.(1996)Fan, Gijbels, Hu and Huang}]{fghh96}
\bibinfo{author}{J.~Fan}, \bibinfo{author}{I.~Gijbels}, \bibinfo{author}{T.C.
  Hu}, \bibinfo{author}{L.S. Huang}, \bibinfo{title}{A study of variable
  bandwidth selection for local polynomial regression},
  \bibinfo{journal}{Statist. Sinica} \bibinfo{volume}{6} (\bibinfo{year}{1996})
  \bibinfo{pages}{113--127}.
\bibitem[{Ferreira et~al.(1997)Ferreira, N{\'u}{\~n}ez-Ant{\'o}n and
  Rodr{\'{\i}}guez-P{\'o}o}]{fnr97}
\bibinfo{author}{E.~Ferreira}, \bibinfo{author}{V.~N{\'u}{\~n}ez-Ant{\'o}n},
  \bibinfo{author}{J.~Rodr{\'{\i}}guez-P{\'o}o}, \bibinfo{title}{Kernel
  regression estimates of growth curves using nonstationary correlated errors},
  \bibinfo{journal}{Statist. Probab. Lett.} \bibinfo{volume}{34}
  (\bibinfo{year}{1997}) \bibinfo{pages}{413--423}.
\bibitem[{Francisco-Fern{\'a}ndez et~al.(2004)Francisco-Fern{\'a}ndez, Opsomer
  and Vilar-Fern{\'a}ndez}]{fov04}
\bibinfo{author}{M.~Francisco-Fern{\'a}ndez}, \bibinfo{author}{J.~Opsomer},
  \bibinfo{author}{J.M. Vilar-Fern{\'a}ndez}, \bibinfo{title}{Plug-in bandwidth
  selector for local polynomial regression estimator with correlated errors},
  \bibinfo{journal}{J. Nonparametr. Stat.} \bibinfo{volume}{16}
  (\bibinfo{year}{2004}) \bibinfo{pages}{127--151}.
\bibitem[{Francisco-Fern{\'a}ndez and Vilar-Fern{\'a}ndez(2001)}]{fv01}
\bibinfo{author}{M.~Francisco-Fern{\'a}ndez}, \bibinfo{author}{J.M.
  Vilar-Fern{\'a}ndez}, \bibinfo{title}{Local polynomial regression estimation
  with correlated errors}, \bibinfo{journal}{Comm. Statist. Theory Methods}
  \bibinfo{volume}{30} (\bibinfo{year}{2001}) \bibinfo{pages}{1271--1293}.
\bibitem[{Gasser and M{\"u}ller(1984)}]{gm84}
\bibinfo{author}{T.~Gasser}, \bibinfo{author}{H.G. M{\"u}ller},
  \bibinfo{title}{Estimating regression functions and their derivatives by the
  kernel method}, \bibinfo{journal}{Scand. J. Statist.} \bibinfo{volume}{11}
  (\bibinfo{year}{1984}) \bibinfo{pages}{171--185}.
\bibitem[{Hall et~al.(1995)Hall, Lahiri and Polzehl}]{hlp95}
\bibinfo{author}{P.~Hall}, \bibinfo{author}{S.N. Lahiri},
  \bibinfo{author}{J.~Polzehl}, \bibinfo{title}{On bandwidth choice in
  nonparametric regression with both short- and long-range dependent errors},
  \bibinfo{journal}{Ann. Statist.} \bibinfo{volume}{23} (\bibinfo{year}{1995})
  \bibinfo{pages}{1921--1936}.
\bibitem[{Hart and Wehrly(1986)}]{hw86}
\bibinfo{author}{J.D. Hart}, \bibinfo{author}{T.E. Wehrly},
  \bibinfo{title}{Kernel regression estimation using repeated measurements
  data}, \bibinfo{journal}{J. Amer. Statist. Assoc.} \bibinfo{volume}{81}
  (\bibinfo{year}{1986}) \bibinfo{pages}{1080--1088}.
\bibitem[{Hart and Wehrly(1993)}]{hw93}
\bibinfo{author}{J.D. Hart}, \bibinfo{author}{T.E. Wehrly},
  \bibinfo{title}{Consistency of cross-validation when the data are curves},
  \bibinfo{journal}{Stochastic Process. Appl.} \bibinfo{volume}{45}
  (\bibinfo{year}{1993}) \bibinfo{pages}{351--361}.
\bibitem[{Masry(2003)}]{m03}
\bibinfo{author}{E.~Masry}, \bibinfo{title}{Local polynomial fitting under
  association}, \bibinfo{journal}{J. Multivariate Anal.} \bibinfo{volume}{86}
  (\bibinfo{year}{2003}) \bibinfo{pages}{330--359}.
\bibitem[{Masry and Fan(1997)}]{mf97}
\bibinfo{author}{E.~Masry}, \bibinfo{author}{J.~Fan}, \bibinfo{title}{Local
  polynomial estimation of regression functions for mixing processes},
  \bibinfo{journal}{Scand. J. Statist.} \bibinfo{volume}{24}
  (\bibinfo{year}{1997}) \bibinfo{pages}{165--179}.
\bibitem[{Opsomer et~al.(2001)Opsomer, Wang and Yang}]{owy01}
\bibinfo{author}{J.~Opsomer}, \bibinfo{author}{Y.~Wang},
  \bibinfo{author}{Y.~Yang}, \bibinfo{title}{Nonparametric regression with
  correlated errors}, \bibinfo{journal}{Statist. Sci.} \bibinfo{volume}{16}
  (\bibinfo{year}{2001}) \bibinfo{pages}{134--153}.
\bibitem[{P{\'e}rez-Gonz{\'a}lez et~al.(2009)P{\'e}rez-Gonz{\'a}lez,
  Vilar-Fern{\'a}ndez and Gonz{\'a}lez-Manteiga}]{P09}
\bibinfo{author}{A.~P{\'e}rez-Gonz{\'a}lez}, \bibinfo{author}{J.M.
  Vilar-Fern{\'a}ndez}, \bibinfo{author}{W.~Gonz{\'a}lez-Manteiga},
  \bibinfo{title}{Asymptotic properties of local polynomial regression with
  missing data and correlated errors}, \bibinfo{journal}{Ann. Inst. Statist.
  Math.} \bibinfo{volume}{61} (\bibinfo{year}{2009}) \bibinfo{pages}{85--109}.
\bibitem[{Perrin(1999)}]{P99}
\bibinfo{author}{O.~Perrin}, \bibinfo{title}{Quadratic variation for {G}aussian
  processes and application to time deformation}, \bibinfo{journal}{Stochastic
  Process. Appl.} \bibinfo{volume}{82} (\bibinfo{year}{1999})
  \bibinfo{pages}{293--305}.
\bibitem[{Rice and Silverman(1991)}]{rs91}
\bibinfo{author}{J.A. Rice}, \bibinfo{author}{B.W. Silverman},
  \bibinfo{title}{Estimating the mean and covariance structure
  nonparametrically when the data are curves}, \bibinfo{journal}{J. Roy.
  Statist. Soc. Ser. B} \bibinfo{volume}{53} (\bibinfo{year}{1991})
  \bibinfo{pages}{233--243}.
\bibitem[{Ruppert(1997)}]{r97}
\bibinfo{author}{D.~Ruppert}, \bibinfo{title}{Empirical-bias bandwidths for
  local polynomial nonparametric regression and density estimation},
  \bibinfo{journal}{J. Amer. Statist. Assoc.} \bibinfo{volume}{92}
  (\bibinfo{year}{1997}) \bibinfo{pages}{1049--1062}.
\bibitem[{Ruppert et~al.(1995)Ruppert, Sheather and Wand}]{rsw95}
\bibinfo{author}{D.~Ruppert}, \bibinfo{author}{S.J. Sheather},
  \bibinfo{author}{M.P. Wand}, \bibinfo{title}{An effective bandwidth selector
  for local least squares regression}, \bibinfo{journal}{J. Amer. Statist.
  Assoc.} \bibinfo{volume}{90} (\bibinfo{year}{1995})
  \bibinfo{pages}{1257--1270}.
\bibitem[{Ruppert and Wand(1994)}]{rw94}
\bibinfo{author}{D.~Ruppert}, \bibinfo{author}{M.P. Wand},
  \bibinfo{title}{Multivariate locally weighted least squares regression},
  \bibinfo{journal}{Ann. Statist.} \bibinfo{volume}{22} (\bibinfo{year}{1994})
  \bibinfo{pages}{1346--1370}.
\bibitem[{Wand and Jones(1995)}]{wj95}
\bibinfo{author}{M.P. Wand}, \bibinfo{author}{M.C. Jones},
  \bibinfo{title}{Kernel smoothing}, volume~\bibinfo{volume}{60} of
  \textit{\bibinfo{series}{Monographs on Statistics and Applied Probability}},
  \bibinfo{publisher}{Chapman and Hall Ltd.}, \bibinfo{address}{London},
  \bibinfo{year}{1995}.
\bibitem[{Yao(2007)}]{y07}
\bibinfo{author}{F.~Yao}, \bibinfo{title}{Asymptotic distributions of
  nonparametric regression estimators for longitudinal or functional data},
  \bibinfo{journal}{J. Multivariate Anal.} \bibinfo{volume}{98}
  (\bibinfo{year}{2007}) \bibinfo{pages}{40--56}.

\end{thebibliography}

\end{document}